\documentclass{lms}

\usepackage{amsmath}
\usepackage{amsfonts}
\usepackage{graphicx,verbatim,color}
\newcommand{\R}{\mathbb{R}}

\newcommand{\Z}{\mathbb{Z}}
\newcommand{\N}{\mathbb{N}}
\newcommand{\hd}{\dim_{\textup{H}}}
\newcommand{\bd}{\dim_{\textup{B}}}
\newcommand{\ubd}{\overline{\dim}_{\textup{B}}}
\newcommand{\lbd}{\underline{\dim}_{\textup{B}}}

\newcommand{\si}{S_{\mathbf{i}}}

\newcommand{\I}{\mathcal{I}}

\newtheorem{theorem}{Theorem}[section] 
\newtheorem{lemma}[theorem]{Lemma}     
\newtheorem{corollary}[theorem]{Corollary}
\newtheorem{proposition}[theorem]{Proposition}
\newnumbered{assertion}{Assertion}    
\newnumbered{conjecture}{Conjecture} 
\newnumbered{definition}{Definition}
\newnumbered{hypothesis}{Hypothesis}
\newnumbered{remark}{Remark}
\newnumbered{note}{Note}
\newnumbered{observation}{Observation}
\newnumbered{problem}{Problem}
\newnumbered{question}{Question}
\newnumbered{algorithm}{Algorithm}
\newnumbered{example}{Example}
\newunnumbered{notation}{Notation}
\title[Inhomogeneous self-affine sets]{The dimensions of inhomogeneous self-affine sets} 
\author{Stuart A. Burrell \and Jonathan M. Fraser}
\dedication{}
\extraline{}
\allowdisplaybreaks
\begin{document}
\maketitle

\begin{abstract}
We prove that the upper box dimension of an inhomogeneous self-affine set is bounded above by the maximum of the affinity dimension and the dimension of the condensation set. In addition, we determine sufficient conditions for this upper bound to be attained, which, in part, constitutes an exploration of the capacity for the condensation set to mitigate dimension drop between the affinity dimension and the corresponding homogeneous attractor. Our work improves and unifies  previous results on general inhomogeneous attractors,  low-dimensional affine systems, and inhomogeneous self-affine carpets, while providing inhomogeneous analogues of Falconer's seminal results on homogeneous self-affine sets. \\


\end{abstract}

\section{Introduction}\label{Introduction} 
A map $S : \R^n \rightarrow \R^n$ is \emph{affine} if it can be written
$$
S(x) = Ax + b
$$
for $A \in \textnormal{GL}(\R, n)$ and translation vector $b \in \R^n$, and is \emph{contracting} if there exists  $c \in (0, 1)$ such that 
$$
|S(x) - S(y)| \leq c|x - y|
$$
for all $x, y \in \R^n$. An affine  \emph{iterated function system} (IFS) is a finite collection $\{S_i\}_{i = 1}^{N}$ of contracting affine maps.  A classic application of Banach's contraction mapping theorem (for details, see \cite{Falconer}) proves the existence of a unique non-empty compact set $F$, called a \emph{homogeneous} attractor, or self-affine set, such that
$$
F = \bigcup\limits_{i = 1}^{N}S_i(F).
$$
There is a natural generalisation of this construction.  If we fix a compact set $C \subset \mathbb{R}^n$, then there exists a unique non-empty compact set $F_C$ such that
$$
F_C = \bigcup\limits_{i = 1}^{N}S_i(F_C) \cup C,
$$
called the \emph{inhomogeneous} attractor, or inhomogeneous self-affine set, with \emph{condensation set} $C$. To express $F_C$ in an amenable way, we require some notation. Henceforth, let $\mathbb{I} = \{S_i\}_{i = 1}^{N}$ denote an affine IFS and $\mathcal{I} = \{1,\dots,N\}$. We write $\si = S_{i_1} \circ \dots \circ S_{i_k}$ for $\textbf{i} = (i_1,i_2,\dots,i_k) \in \mathcal{I}^k$. Furthermore, let
$$
\mathcal{I}^* = \bigcup\limits_{k = 1}^{\infty} \mathcal{I}^k
$$
denote the set of finite words over $\mathcal{I}$. An elegant description of $F_C$, from \cite{Fraser} and \cite{sniphd}, is
$$
F_C = F_{\emptyset} \cup \mathcal{O},
$$
where  $F_\emptyset$ is the homogeneous attractor (corresponding to $C = \emptyset$), and $\mathcal{O}$ is the \emph{orbital set} defined by
$$
\mathcal{O} = C \cup \bigcup\limits_{\textbf{i} \in \mathcal{I}^*} \si(C).
$$
 Since their introduction by Barnsley (1985) \cite{barn} and Hata (1985) \cite{hata}, inhomogeneous attractors have received further attention in, for example,  \cite{fraserbaker,burrell,Fraser,18,antti,olssni,sniphd}. A natural question explored in recent work concerns the relationship between the dimensions of $F_C$, $C$ and $F_\emptyset$. In particular, one may wonder in what situations
\begin{equation}\label{maineq}
\dim F_C = \max\left\{\dim F_\emptyset, \dim C\right\},
\end{equation}
where $\dim$ denotes some notion of dimension. For dimensions satisfying countable stability, such as the Hausdorff or packing dimension, this is immediate. Consequently, the recent focus has been on box dimension, a popular example of a countably unstable dimension. Recall that for a non-empty bounded set $F\subseteq \R^n$, the upper and lower box dimensions are defined as
$$
\ubd F = \limsup_{\delta \rightarrow 0} \frac{\log N_{\delta}(F)}{-\log \delta},
$$
and
$$
\lbd F = \liminf_{\delta \rightarrow 0} \frac{\log N_{\delta}(F)}{-\log \delta},
$$
respectively, where $N_{\delta}(F)$ denotes the minimum number of hypercubes of sidelength $\delta$ required to cover $F$. If these values coincide, we say the set has box dimension equal to the common value and denote this by $\bd F$.\\

In \cite{fraserbaker,Fraser,olssni,sniphd}, various solutions to (\ref{maineq}) are given for upper box dimension in the case where $\mathbb{I}$ consists of similarity mappings. For systems containing arbitrary bi-Lipschitz maps, bounds on $\ubd F_C$ are given by Burrell based on upper Lipschitz dimension \cite{burrell}. Corollaries of this result establish (\ref{maineq}) for some low-dimensional affine systems and those satisfying bounded distortion, such as conformal systems (see \cite{feng} for definitions). For upper box dimension, (\ref{maineq}) may fail for self-similar sets with overlaps \cite{fraserbaker} and specific self-affine settings \cite{18}. In the case of lower box dimension, \eqref{maineq} fails to hold generally even for self-similar systems satisfying the strong separation condition \cite{Fraser}.\\

The typical strategy used to approach (\ref{maineq}), introduced in \cite{Fraser}, is to establish bounds of the form
\begin{equation}\label{mainform}
\max\left\{\ubd F_\emptyset, \ubd C\right\} \leq \ubd F_C \leq \max\left\{s, \ubd C\right\},
\end{equation}
where $s \in \R$ is  a natural estimate for $\ubd F_\emptyset$, such as similarity dimension in the self-similar case \cite{Fraser} or upper Lipschitz dimension \cite{burrell} in the general case. This exploits the abundance of literature on the equality of $s$ and $\ubd F_\emptyset$ in different settings, which may then determine precise conditions for equality. In  the affine setting, the natural candidate for $s$ is the \emph{affinity dimension}, and we prove (\ref{mainform}) holds in this case.\\

The affinity dimension is derived from the notion of Falconer's singular value function, introduced in \cite{1988kjf}. The singular values of $A \in \textnormal{GL}(\R, n)$ are written $\alpha_j(A)$ (or simply $\alpha_j$) and correspond to the lengths of the mutually perpendicular principal axes of $A(B)$, where $B$ denotes a ball of unit diameter in $\R^n$ \cite{1988kjf}. Alternatively, they are the positive square roots of the eigenvalues of $AA^T$. We adopt the convention $1 > \alpha_1 \geq \alpha_2 \geq \dots \geq \alpha_n > 0$. For $0 \leq s \leq n$, the singular value function of $A \in \textnormal{GL}(\R, n)$ is given by
$$
\phi^s(A) = \alpha_1(A)\alpha_2(A)\cdots\alpha_m(A)^{s-m+1},
$$
where $m \in \Z$ satisfies $m-1 < s \leq m$. Finally, as in \cite{1988kjf}, we  define $\phi^s(A) = (\det A)^{s/n}$ for $s > n$. Moreover, for convenience we set $\phi^s(S) = \phi^s(A)$ where $A$ is the linear component of a general affine map $S$. Then, for each $k\in \N$, define $s_k$ to be the solution of
$$
\sum\limits_{\textbf{i} \in \mathcal{I}^k} \phi^{s_k}(\si) = 1.
$$ 
The corresponding limit, denoted throughout by $s$,
$$
s := \lim_{k \rightarrow \infty} s_k,
$$
exists and is known as the \emph{affinity dimension} associated with $\mathbb{I}$. It is proven in \cite{burrell} that if the affinity dimension $s$ is less than or equal to $1$ and coincides with $\ubd F_\emptyset$, then
\begin{equation}\label{affineeq}
\ubd F_C = \max\left\{\ubd{F_\emptyset}, \ubd C\right\}.
\end{equation}
This is an immediate corollary of \cite[Theorem 2.1]{burrell}, arising from the fact that when the affinity dimension is less than one it coincides with the upper Lipschitz dimension. Otherwise, it is elementary to see that the affinity dimension is generally strictly less than the upper Lipschitz dimension. Thus, establishing (\ref{mainform}) for affinity dimension constitutes a natural and strictly improved bound for affine systems in comparison to the universal bound from \cite{burrell}.

\section{Results}

Our main result may be considered an inhomogeneous analogue of Falconer's seminal result on homogeneous self-affine sets \cite{1988kjf}, which establishes $\ubd F_\emptyset \leq s$.
\begin{theorem}\label{mainaffine}
Let $F_C \subset \mathbb{R}^n$ be an inhomogeneous self-affine set with compact condensation set $C \subset \mathbb{R}^n$. We have
$$
\max\left\{\ubd F_{\emptyset},\, \ubd C\right\} \leq \ubd F_C \leq \max\left\{{s, \ubd C}\right\},
$$
where $s$ is the affinity dimension associated with the underlying IFS.
\end{theorem}

The following corollary is immediate.
\begin{corollary}\label{maincor}
Let $F_C \subset \mathbb{R}^n$ be an inhomogeneous self-affine set with compact condensation set $C \subset \mathbb{R}^n$ and let $s$ be the associated affinity dimension.  Then
\begin{enumerate}
\item if $ \ubd F_\emptyset = s$, then $\ubd F_C = \max\left\{{\ubd F_\emptyset, \ubd C}\right\}$,
\item if $ \ubd C \geq s$, then $\ubd F_C =  \ubd C.$
\end{enumerate}
\end{corollary}

Establishing precise conditions for the affinity dimension to coincide with $\ubd F_\emptyset$ is a major open problem in fractal geometry and has been the focus of considerable amounts of work, for example \cite{bhr,1988kjf,hd2,fraserbox,18,fraserkemp,27,jordanjurga}. Therefore there are numerous explicit and non-explicit situations where Corollary \ref{maincor} provides a precise result, and an affirmative  solution to \eqref{maineq} in the self-affine setting.  For example, a well-known result by Falconer \cite{1988kjf} states that $s = \ubd F_\emptyset= \hd F_\emptyset$ almost surely if one randomises the translation vectors associated with the affine maps, provided the linear parts all have norm strictly bounded  above by 1/2, see also \cite{jordanjurga}.  Falconer proved in a subsequent paper that if $F_\emptyset \subset \R^2$ satisfies some separation conditions and contains a connected component not contained in a straight line, then $s = \ubd F_\emptyset$ holds, see  \cite[Corollary 5]{hd2}. A recent breakthrough result of B\'ar\'any, Hochman and Rapaport \cite{bhr} proves $s = \ubd F_\emptyset = \hd F_\emptyset$ in the planar case assuming only strong separation, together with  mild non-compactness and irreducibility assumptions on the linear components of the maps $S_i$.\\

The next result explores the case where $\ubd F_C   > \max\{\ubd F_\emptyset, \ubd C\}$, that is when \eqref{maineq} fails. This is an exploration of  conditions under which $C$ compensates for dimension drop between $s$ and $\ubd F_{\emptyset}$. To state this result, we require the definition of the \emph{condensation open set condition} (COSC), appearing in \cite{antti,olssni,sniphd} and $m$-$\delta$-\emph{stoppings}. Firstly, an IFS satisfies the COSC if there exists an open set $U$ with 
$$C \subset U \setminus \bigcup\limits_{i = 1}^{N} \overline{S_i(U)},$$ such that $S_i(U)\subset U$ for $i = 1,\dots,N$, and
$i \neq j \implies S_i(U) \cap S_j(U) = \emptyset$. Secondly, for each $1 \leq m \leq n$ and $\delta \in (0,1]$, define the $m$-$\delta$-stopping to be
$$
\mathcal{I}_m(\delta) = \{\textbf{i} \in \I^* : \alpha_m(\si) < \delta \leq \alpha_m(S_{\textbf{i}_-})\},
$$
where $\textbf{i}_- = (i_1,\dots,i_{k-1})$ for $\textbf{i} = (i_1,\dots,i_k)$.  For the next theorem we will only use  $\mathcal{I}_n(\delta)$, but later in the proofs section we will use it more generally and so introduced it here in full generality for brevity. Throughout, we fix a compact ball $X \subset \R^n$ such that $S_i(X) \subset X$ for $i = 1, \dots, N$ and $C \subseteq X$. Such a ball always exists and without loss of generality, we may assume that $X$ has unit diameter. 
\begin{theorem}\label{projectivethm}
Let $\mathbb{I} = \{S_i\}_{i = 1}^{N}$ denote an affine IFS with affinity dimension $s \leq n$ and condensation set $C \subset \R^n$ satisfying the COSC. If $\lbd C \geq n-1$ and there exists $\kappa > 0$ such that for all $\delta \in (0, 1]$ and $\textbf{i} \in \I_n(\delta)$ we have
$$
N_{\delta}(\si(C)) \geq \kappa N_{\delta}(\si(X)),
$$
then
$$
\ubd F_C = \max\left\{s, \ubd C\right\}
$$
and
$$
 \max\left\{s, \lbd C \right\} \leq \lbd F_C \leq  \max\left\{s, \ubd C \right\}.
$$
\end{theorem}
Note that the condition of the theorem is independent of the choice of ball $X$, although the constant $\kappa$ may change. The fact that we only get bounds for the lower box dimension of $F_C$ should not come as a surprise and one should not expect to be able to improve these bounds in general, see \cite{Fraser}.  Note that if, in the setting of Theorem \ref{projectivethm}, the box dimension of $C$ exists, then so does the box dimension of $F_C$.\\

The assumption in Theorem \ref{projectivethm} arises in quite natural circumstances, for example, the setting of the following proposition, an inhomogeneous analogue of Falconer's  \cite[Proposition 4]{hd2}, requires only that $C$ be in some sense robust under projection onto subspaces. Let $\mathcal{L}^{k}$ denote $k$-dimensional Lebesgue measure and $P_k$ denote the set of orthogonal projections onto $k$-dimensional subspaces of $\R^n$.
\begin{proposition}\label{projprop}
Let $F_C \subset \mathbb{R}^n$ be an inhomogeneous self-affine set with compact condensation set $C \subset \mathbb{R}^n$ satisfying the COSC and let $s \leq n$ be the associated affinity dimension. If
$$
\inf\limits_{\pi\in P_{n-1}}\mathcal{L}^{n-1}(\pi C) > 0,
$$
then
$$
\ubd F_C = \max\left\{s, \ubd C\right\}
$$
and
$$
 \max\left\{s, \lbd C \right\} \leq \lbd F_C \leq  \max\left\{s, \ubd C \right\}.
$$
\end{proposition}
The robustness  assumption on $C$ in Proposition  \ref{projprop} forces $\lbd C \geq n-1 $ and so this result only yields new information when $s > n-1 $. It is interesting to compare this result with \cite[Corollary 2.3]{burrell} and the discussion thereafter, which applies to self-affine systems where $s \leq 1$. \\

The projection of a connected set in $\mathbb{R}^2$ which is not contained in a line onto a line contains an interval with length uniformly bounded away from 0. This observation yields the following corollary of Proposition \ref{projprop}.
\begin{corollary}\label{projcor}
Let $F_C \subset \R^2$ be an inhomogeneous self-affine set with compact condensation set $C \subset \mathbb{R}^2$ satisfying the COSC and let $s \leq 2$ be the associated affinity dimension. If $C$ has a connected component not contained in a line, then
$$
\ubd F_C = \max\left\{s, \ubd C\right\}
$$
and
$$
 \max\left\{s, \lbd C \right\} \leq \lbd F_C \leq  \max\left\{s, \ubd C \right\}.
$$
\end{corollary}
 The reader may find it interesting to notice the parallels between this result and  Falconer's \cite[Corollary 5]{hd2}, which concerns  the equality of $\ubd F_\emptyset$ and $s$ under similar conditions concerning the robustness of connected components under projection.  In some sense our  inhomogeneous analogue is easier to use than the homogeneous result of Falconer.  Our result requires a connectedness condition on $C$, which is given, whereas the homogeneous result requires one to check a connectedness condition on $F_\emptyset$, which depends delicately on the IFS.  Moreover, the separation assumption makes it difficult for  $F_\emptyset$ to be connected at all.  For example, the  strong separation condition forces $F_\emptyset$ to be totally disconnected, but our result can still apply in this setting.\\

The above results provide new families of inhomogeneous attractors where \eqref{maineq} fails for the upper (and lower) box dimension.  We illustrate this by example.  Let $n=2$ and $\mathbb{I} =\{S_1, S_2\}$, where $S_1, S_2$ are the linear maps associated with the matrices
\[
\begin{bmatrix}
  1/2       & 0 \\
 1/2    &  1/2
\end{bmatrix}, \qquad  \begin{bmatrix}
  1/2       &  1/2 \\
0   &  1/2
\end{bmatrix}
\]
respectively.  It is clear that the affinity dimension of this system is strictly greater than one and that $F_\emptyset$ is just a single point at the origin.  Let $C$ be the boundary of a circle centred at $(3/4, 3/4)$ with radius $1/5$.  It is also clear that the COSC is satisfied by taking $U = (0,1)^2$ and that $C$ is connected but not contained in a line, see Figure 1.  It follows from Corollary \ref{projcor} that
\[
\ubd F_C = \lbd F_C  = s > 1 =  \max\left\{\bd  F_\emptyset, \bd C\right\}.
\]
\begin{figure}[h]\label{ovalsfig}
\includegraphics[scale=0.85]{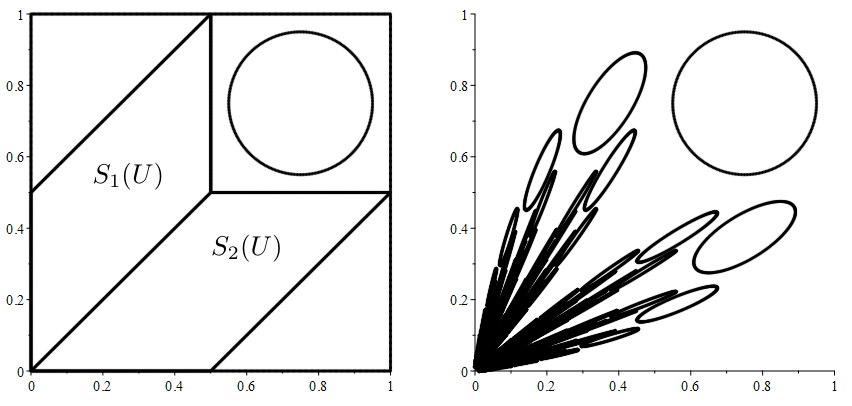}
\caption{A bouquet of ovals: the condensation set together with the two images of the open rectangle $U = (0,1)^2$ (left) and the corresponding inhomogeneous self-affine set (right).}
\end{figure}
This is the first counter example to \eqref{maineq} where  $F_\emptyset$ is a single point and the OSC is satisfied.  Moreover, it was shown in \cite[Corollary 4.9]{fraserbaker} that for planar inhomogeneous self-\emph{similar} sets one always has
\[
\ubd F_C \leq \max\left\{ \ubd C, \ \ubd F_\emptyset+\ubd C-\frac{\ubd F_\emptyset\ubd C}{s} \right\},
\]
where $s$ is the similarity dimension.  In particular this shows that when $\ubd F_\emptyset = 0$ the formula \eqref{maineq} cannot fail.  The example presented above shows that this phenomenon does \emph{not} extend to the self-affine case. It was also shown in  \cite[Corollary 4.8]{fraserbaker} that, in the self-similar setting, if $\max\{\ubd  F_\emptyset, \ubd C\} < s$, then $\ubd F_C < s$.  The above example also demonstrates that this does not extend to the self-affine setting. \\

The assumption in Proposition \ref{projprop} is by no means necessary, and advancements in the homogeneous setting may illuminate further the capacity for $C$ to mitigate dimension drop. Excitingly, we suggest the natural interplay between these questions may
allow further study of inhomogeneous attractors to 
translate into novel conditions relating to dimension drop in the homogeneous case. Specifically, this may arise from solutions to the following.

\begin{question}
Consider an affine IFS $\mathbb{I} = \{S_i\}_{i = 1}^{N}$ with condensation set $C \subset \R^n$. If $s > \ubd F_\emptyset$, then what conditions guarantee
$$
\ubd F_C = \max\left\{s, \ubd C\right\}?
$$
\end{question}

\section{Proof of Theorem \ref{mainaffine}}

Let $\mathbb{I} = \{S_i\}_{i=1}^{N}$ be an affine IFS and $C\subseteq X$ be compact. Denote the affinity dimension of $\mathbb{I}$ by $s$ and assume $s \leq n$, since if $s > n$ the result is trivial. \\

It follows immediately from the definition of box dimension that for $t > \ubd C$ there exists a constant $A_t$ satisfying
\begin{equation}\label{eqN(c)}
N_{\delta}(C) \leq A_t\delta^{-t}
\end{equation} 
for all $\delta \in (0, 1]$. In addition, if $t > s$, then 
\begin{equation}\label{phicnst}
B_t := \sum\limits_{\textbf{i} \in \I^*}\phi^t(\si)  < \infty
\end{equation}
by \cite[Proposition 4.1 (c)]{1988kjf}, where $B_t$ depends only on $t$. We fix a constant $b \in \R$ satisfying
$$
0 < b < \min\limits_{i = 1,\dots,N} \alpha_n(S_i)< 1,
$$
and note for any $\delta \in (0, 1]$, $1 \leq m \leq n$ and $\textbf{i} \in \I_m(\delta)$, we have
\begin{equation}\label{bt}
\delta \geq \alpha_m(\si) \geq \alpha_m(S_{\textbf{i}_-})b \geq \delta b.
\end{equation}

Prior to reading the subsequent arguments, the following simple geometric observation, employed frequently in our proofs, may aid the reader less familiar with the classical arguments on self-affine sets found in \cite{1988kjf} or \cite{Falconer}. Consider an ellipsoid $E$ with principal axes of lengths $l_1, \dots, l_n$. For dimension calculations, we are interested in obtaining an estimate of the number of hypercubes of a given sidelength required to cover such ellipsoids. Constants are typically inconsequential, and so often a coarse estimate suffices. The minimum number of hypercubes of sidelength $l_m$ required to cover $E$ is at most
\begin{equation}\label{coverbound}
\left(\frac{l_1}{l_m}+1\right)\left(\frac{l_2}{l_m} + 1\right)\cdots\left(\frac{l_{m-1}}{l_m}+1\right) \leq 2^n \frac{l_1}{l_m}\frac{l_2}{l_m}\cdots\frac{l_{m-1}}{l_m} = 2^nl_1l_2 \cdots l_{m-1}l_m^{-m + 1}.
\end{equation}
This can be seen by first covering $E$ by a minimal hypercuboid of sidelengths equal to the principal axes of $E$ and then covering this optimally. Figure 2 illustrates this fact for a cuboid of sidelengths $a > b > c$ in $\R^3$. Specifically, we see that $2a/b$ cubes of sidelength $b$ would suffice, whereas we would require a single cube of sidelength $a$ or at most $2^2(a/c)(b/c)$ cubes of sidelength $c$.
\begin{figure}[h]\label{coverpic}
\begin{center}
\includegraphics[scale = 0.25]{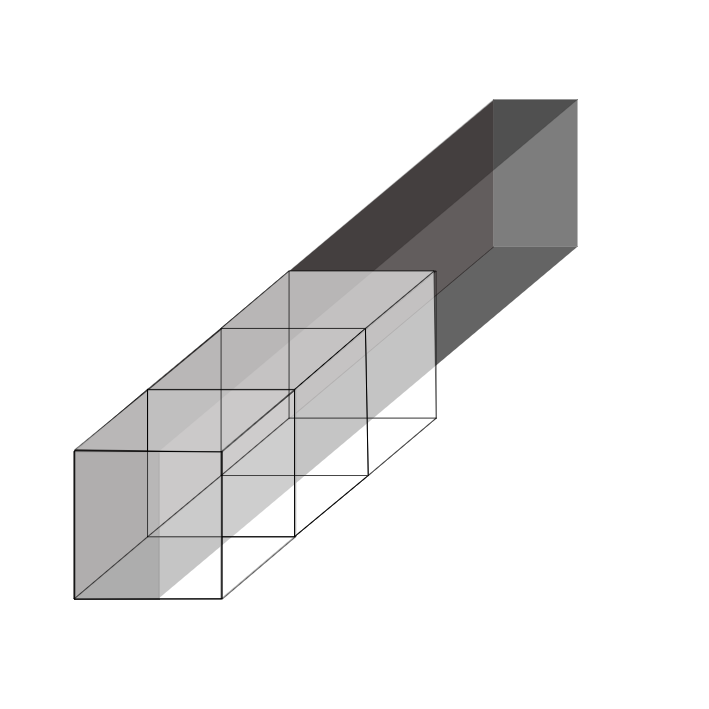}
\caption{Covering a cuboid of sidelengths $a > b > c$ in $\R^3$ with cubes of sidelength $b$.}
\end{center}
\end{figure}

\subsection{Preliminary lemmas}
\begin{lemma}\label{inclusionlemma}
For $\delta \in (0, 1]$ and $1 \leq m \leq n$, we have
$$
\bigcup\limits_{\substack{\mathbf{i} \in \mathcal{I}^* \\ \delta > \alpha_m(S_{\mathbf{i}})}} S_{\mathbf{i}}(C) \subseteq  \bigcup\limits_{\mathbf{i} \in \mathcal{I}_m(\delta)}S_{\mathbf{i}}(X).
$$
\begin{proof}
For $$x \in \bigcup\limits_{\substack{\textbf{i} \in \mathcal{I}^* \\ \delta > \alpha_m(S_{\textbf{i}})}} S_{\textbf{i}}(C),$$ there exists some $\textbf{i} = (i_1,\dots,i_k) \in \I^*$ such that $x \in \si(C)$ and $\delta > \alpha_m(\si)$. Since $\delta > \alpha_m(\si)$, there also exists some prefix $\textbf{i}_p$ of $\textbf{i}$ with $\textbf{i}_p \in \mathcal{I}_m(\delta)$, and so let us consider the concatenation $\textbf{i} = \textbf{i}_p\textbf{j}$. If $\textbf{j} = \emptyset$, then $\textbf{i} \in \mathcal{I}_m(\delta)$. Else, there exists some $c$ such that $x = \si(c) = S_{\textbf{i}_p}(S_{\textbf{j}}(c)) \in S_{\textbf{i}_p}(X)$ as required.
\end{proof}
\end{lemma}

\begin{lemma}\label{alphasmall}
Fix $1 \leq m \leq n$ and let $\mathbf{i} \in \I^*$ be such that $\alpha_m(\si) < \delta$. We have
$$
N_{\delta}(\si(X)) \leq 2^{n} \frac{\alpha_1(\si)}{\alpha_m(\si)}\frac{\alpha_2(\si)}{\alpha_m(\si)}\cdots\frac{\alpha_{m-1}(\si)}{\alpha_m(\si)}.
$$
\begin{proof}
First note that $\si(X)$ is  an ellipsoid with principal  axes having lengths equal to the singular values of $\si$. The result then follows follows immediately from the geometric observation described by equation (\ref{coverbound}).
\end{proof}
\end{lemma}

\begin{lemma}\label{alphabig}
Let $\ubd C \leq t \leq n$ and $m \in \Z$ be the integer satisfying $m-1 < t \leq m$. If $\mathbf{i} \in \I^*$ is such that $\alpha_m(\si) \geq \delta$, then
$$
N_{\delta}(\si(C)) \leq 2^nA_t\delta^{-t}\phi^t(\si),
$$
where $A_t$ is a constant depending only on $t$.
\begin{proof}
The image under $\si$ of a cover of $C$ by balls of diameter $\delta/\alpha_m(\si)$ is a cover of $\si(C)$ by ellipsoids with the $m$ largest principal axes of lengths 
$$
\alpha_i(\si)\left(\frac{\delta}{\alpha_m(\si)}\right) = \delta\frac{\alpha_i(\si)}{\alpha_m(\si)} 
$$
for $i = 1,\dots, m$, the smallest of which has length $\delta$. Each such ellipsoid can be covered by at most
$$
\frac{2\delta\frac{\alpha_1(\si)}{\alpha_m(\si)}}{\delta}\frac{2\delta\frac{\alpha_2(\si)}{\alpha_m(\si)}}{\delta} \cdots \frac{2\delta\frac{\alpha_{m-1}(\si)}{\alpha_m(\si)}}{\delta}  \leq 2^n\frac{\alpha_1(\si)}{\alpha_m(\si)}\frac{\alpha_2(\si)}{\alpha_m(\si)}\cdots\frac{\alpha_{m-1}(\si)}{\alpha_m(\si)}
$$
hypercubes of sidelength $\delta$. Hence
\begin{align*}
N_{\delta}(\si(C)) &\leq N_{\delta/\alpha_m(\si)}(C)\left( 2^n\frac{\alpha_1(\si)}{\alpha_m(\si)}\frac{\alpha_2(\si)}{\alpha_m(\si)}\cdots\frac{\alpha_{m-1}(\si)}{\alpha_m(\si)}\right)\\ &\leq A_t \left(\frac{\delta}{\alpha_m(\si)}\right)^{-t}\left(2^n\frac{\alpha_1(\si)}{\alpha_m(\si)}\frac{\alpha_2(\si)}{\alpha_m(\si)}\cdots\frac{\alpha_{m-1}(\si)}{\alpha_m(\si)}\right)\\ &= 2^nA_t\delta^{-t}\phi^t(\si)
\end{align*}
as required.
\end{proof}
\end{lemma}
\subsection{Proof of Theorem \ref{mainaffine}}

Monotonicity and finite stability of upper box dimension imply 
$$
\max\left\{\ubd F_\emptyset, \ubd C\right\} \leq \ubd F_C \leq \max\left\{\ubd F_\emptyset, \ubd \mathcal{O}\right\}
$$
and so it suffices to show that
$$
\ubd \mathcal{O} \leq \max\left\{s, \ubd C\right\}
$$
since it is well known (see \cite[Theorem 9.12]{Falconer}) that $s \geq \ubd F_\emptyset$. Fix $\delta \in (0, 1]$ and $t > \max\{s, \ubd C\}$. If $\max\{s, \ubd C\} \geq n$ then the result is trivial, so we may assume $t \leq n$. For $m \in \Z$ satisfying $m - 1 < t \leq m$, we have
\begin{flalign*}
\delta^t N_{\delta}(\mathcal{O}) &= \delta^tN_{\delta}\left(C \cup \bigcup\limits_{\textbf{i} \in \mathcal{I}^*}  S_{\textbf{i}}(C)\right) \\
%%%%%%%%%%%%%%%%%%%%%%%%%%%%%%%%%%%%%%%%%%%%%%%%%%%%%%
%& \shortintertext{[Since $\mathcal{O} = C \cup \bigcup\limits_{\textbf{i} \in \mathcal{I}^*}  S_{\textbf{i}}(C)$.]}
& \leq A_t + \delta^tN_{\delta}\left(\bigcup\limits_{\substack{\textbf{i} \in \mathcal{I}^* \\\alpha_m(S_{\textbf{i}}) \geq \delta}} S_{\textbf{i}}(C)\right) + 
\delta^tN_{\delta}\left(\bigcup\limits_{\substack{\textbf{i} \in \mathcal{I}^* \\\alpha_m(S_{\textbf{i}}) < \delta}} S_{\textbf{i}}(C)\right) \textnormal{\,\,\,\,(using (\ref{eqN(c)}))}\\
& \leq A_t + \delta^t\sum\limits_{\substack{\textbf{i} \in \mathcal{I}^* \\\alpha_m(S_{\textbf{i}}) \geq \delta}} N_{\delta}(S_{\textbf{i}}(C)) + 
\delta^t\sum\limits_{\textbf{i} \in \mathcal{I}_m(\delta)} N_{\delta}\left(S_{\textbf{i}}(X)\right) \textnormal{\,\,\,\,(by Lemma \ref{inclusionlemma})}\\
&\leq A_t + \delta^t\sum\limits_{\substack{\textbf{i} \in \mathcal{I}^* \\ \alpha_m(S_{\textbf{i}}) \geq \delta}} 2^nA_t\delta^{-t}\phi^t(\si) \\&\,\,\,\,\,\,\,+ 
\delta^t\sum\limits_{\textbf{i} \in \mathcal{I}_m(\delta)} 2^n\frac{\alpha_1(\si)}{\alpha_m(\si)}\frac{\alpha_2(\si)}{\alpha_m(\si)}\cdots\frac{\alpha_{m-1}(\si)}{\alpha_m(\si)} \textnormal{\,\,\,\,(by Lemmas \ref{alphasmall} and \ref{alphabig})}\\
&\leq A_t + 2^nA_t\sum\limits_{\substack{\textbf{i} \in \mathcal{I}^* \\ \alpha_m(S_{\textbf{i}}) \geq \delta}} \phi^t(\si) \\&\,\,\,\,\,\,\,+
2^n\sum\limits_{\textbf{i} \in \mathcal{I}_m(\delta)} \frac{\alpha_1(\si)}{\alpha_m(\si)}\frac{\alpha_2(\si)}{\alpha_m(\si)}\cdots\frac{\alpha_{m-1}(\si)}{\alpha_m(\si)}\frac{\alpha_m(\si)^t}{b^t}\textnormal{\,\,\,\,(using (\ref{bt}))} \\
&\leq A_t + 2^nA_t\sum\limits_{\substack{\textbf{i} \in \mathcal{I}^* \\ \alpha_m(S_{\textbf{i}}) \geq \delta}} \phi^t(\si) + 
\frac{2^n}{b^t}\sum\limits_{\textbf{i} \in \mathcal{I}_m(\delta)} \phi^t(\si) \\
&\leq A_t + 2^nB_t\left(A_t+ b^{-t}\right) \textnormal{\,\,\,\,(using (\ref{phicnst}))}
\end{flalign*}
Thus, 
$$
\frac{\log N_\delta(\mathcal{O})}{-\log \delta} \leq t + \frac{\log\left(A_t + 2^nB_t\left(A_t+ b^{-t}\right)\right)}{-\log \delta},
$$
from which the result follows as $\delta \rightarrow 0$. \hfill $\square$

\section{Proof of Theorem \ref{projectivethm}}
 
Fix $\delta \in (0,1)$ and recall that $s$ denotes the affinity dimension of $\mathbb{I}$. It is stated in \cite{hd2} that for $t < s$ there exists $c_t > 0$ with
\begin{equation}\label{lowerbound}
\sum\limits_{\I_n(\delta)} \phi^t(\si) \geq c_t
\end{equation}
for some constant $c_t$ that does not depend on $\delta$. This follows immediately from \cite[Proposition 4.1 (a)]{1988kjf}. Since we assume $\lbd C \geq n - 1$, if $s \leq n - 1$, then Theorem \ref{mainaffine} implies that $\ubd F_C = \ubd C = \max\{s, \ubd C\}$, and also $\lbd F_C \geq \lbd C = \max\{s, \lbd C\}$.  Thus, henceforth we assume that $n-1 <t< s \leq n$.\\

Let $U$ denote the open set satisfying the COSC. Compactness of $C$ implies that there exists some constant $\eta > 0$ with 
$$
\inf \left\{|x - y| : x\in C,\,\,y \in \bigcup\limits_{i = 1}^{N} \si(U) \cup (\R^n \setminus U)\right\} = 2\eta.
$$
Let $B(C, \eta)$ denote a closed $\eta$-neighborhood of $C$ and $E$ be a hypercube of sidelength $\delta$ in a minimal $\delta$-cover of $\mathcal{O}$. For $\textbf{i} \in \mathcal{I}_n(\delta)$, we have $\si(B(C, \eta))$ is a neighborhood of $\si(C)$ satisfying
$$
\si(B(C, \eta)) \cap F_C = \si(C)
$$
and 
$$
\inf \{|x - y| : x \in \si(C), y\notin \si(B(C, \eta)) \} \geq \alpha_n(\si)\eta > b\delta\eta,
$$
implying
$$
\inf \{|x - y| : x \in \si(C), y\in S_{\textbf{j}}(C) \textnormal{\,\,such that\,\,} \textbf{i},\textbf{j}\in \I_n(\delta), \textbf{i} \neq \textbf{j}\} > 2b\delta\eta.
$$
Let $V_n$ denote the constant such that the volume of an $n$-sphere of radius $2b\eta\delta$ is $V_n\delta^n$. For the sets in $\{\si (C): \textbf{i} \in \mathcal{I}_n(\delta)\}$ that intersect $E$ we can associate pairwise disjoint open sets in $E$ of volume at least $V_n\delta^n/2^n$ (with this lower bound obtained at the vertices) and it therefore follows by a simple volume argument that  $E$ can intersect at most 
$$
\frac{\delta^n}{\frac{1}{2^n}V_n\delta^n} = (2^{-n}V_n)^{-1}
$$
of the sets $\{\si (C): \textbf{i} \in \mathcal{I}_n(\delta)\}$. Hence
\begin{equation}\label{sep}
N_{\delta}(\mathcal{O}) \geq 2^{-n}V_n\sum\limits_{\textbf{i} \in \I_n(\delta)} N_\delta(\si(C)). 
\end{equation}
Our assumption on $C$ implies that for $\textbf{i} \in \I_n(\delta)$ we have
\begin{align}\label{count}
N_\delta(\si(C)) &\geq \kappa N_\delta(\si(X)) \nonumber \\
&\geq \kappa b^nN_{b\delta}(\si(X)) \nonumber\\
&\geq \kappa b^nN_{\alpha_n(\si)}(\si(X))\nonumber\\
&\geq \kappa b^nc\frac{\alpha_1(\si)}{\alpha_n(\si)}\frac{\alpha_2(\si)}{\alpha_n(\si)}\cdots\frac{\alpha_{n-1}(\si)}{\alpha_n(\si)}
\end{align}
for some constant $c>0$ only depending on $n$. This yields
\begin{align*}
N_{\delta}(\mathcal{O}) &\geq 2^{-n}V_n\sum\limits_{\textbf{i} \in \I_n(\delta)} N_\delta(\si(C))  \textnormal{\,\,\,\,(using (\ref{sep}))} \\
&\geq 2^{-n}V_n\sum\limits_{\textbf{i} \in \I_n(\delta)} \kappa b^n c \frac{\alpha_1(\si)}{\alpha_n(\si)}\frac{\alpha_2(\si)}{\alpha_n(\si)}\cdots\frac{\alpha_{n-1}(\si)}{\alpha_n(\si)} \textnormal{\,\,\,\,(using (\ref{count}))}\\
&=\kappa b^n c 2^{-n}V_n\sum\limits_{\textbf{i} \in \I_n(\delta)} \phi^t(\si)\alpha_n(\si)^{-t} \\
&\geq \kappa  b^n c2^{-n}V_n\delta^{-t} \sum\limits_{\textbf{i} \in \I_n(\delta)} \phi^t(\si) \\
&\geq \kappa  b^n c 2^{-n}V_n c_t\delta^{-t} \textnormal{\,\,\,\,(by (\ref{lowerbound}))}.
\end{align*}
Hence $\lbd \mathcal{O} \geq t$, from which it follows that $\ubd F_C \geq \lbd F_C \geq \lbd \mathcal{O} \geq s$, proving the theorem. \hfill $\square$

\section{Proof of Proposition \ref{projprop}}
Let $\mathbb{I} = \{S_i\}_{i = 1}^{N}$ denote an affine IFS with compact condensation set $C \subseteq \R^n$ satisfying the COSC. Moreover, suppose
$$
 \inf\limits_{\pi \in P_{n-1}} \mathcal{L}^{n-1}(\pi C) > 0.
$$
By Theorem \ref{projectivethm} it suffices to show that there exists $\kappa > 0$ such that for all $\delta > 0$ and $\textbf{i} \in \I_n(\delta)$ we have
$$
N_{\delta}(\si(C)) \geq \kappa N_\delta(\si(X)) .
$$
Therefore, in order to reach a contradiction, assume that for arbitrarily small $\kappa>0$ we can find $\delta > 0$ and $\textbf{i} \in \I_n(\delta)$ such that
$$
N_{\delta}(\si(C)) < \kappa N_\delta(\si(X)) \leq \kappa  2^{n} \frac{\alpha_1(\si)}{\alpha_n(\si)}\frac{\alpha_2(\si)}{\alpha_n(\si)}\cdots\frac{\alpha_{n-1}(\si)}{\alpha_n(\si)},
$$
where the final inequality comes from Lemma \ref{alphasmall}. Let $\{E_j\}_j$ be an optimal  cover of $\si(C)$ by hypercubes of sidelength $\delta$ and place each $E_j$ inside a ball $B_j$ of diameter $\sqrt{n} \delta$ and consider  $\{\si^{-1}B_j\}_j$ which is a cover of $C$ by ellipsoids with axes of length   $\sqrt{n}\delta/\alpha_1(\si), \dots, \sqrt{n}\delta/\alpha_n(\si)$.  Note that the longest axes of each of these ellipsoids are all parallel (by the singular value decomposition theorem, for example) and let $\pi$ denote projection onto the $(n-1)$-dimensional hyperplane orthogonal to the common direction of the longest axes of the ellipsoids $\{\si^{-1}B_j\}_j$.  It follows that  $\{\pi\si^{-1}B_j\}_j$  is a cover of $\pi(C)$ by sets, each of which is easily seen to have $(n-1)$-volume at most
\[
 n^{(n-1)/2} \frac{\delta}{\alpha_1(\si)}\frac{\delta}{\alpha_2(\si)}\cdots\frac{\delta}{\alpha_{n-1}(\si)}
\]
and therefore we can bound the $(n-1)$-volume of $\pi(C)$ above by
\[
\kappa  2^{n} \frac{\alpha_1(\si)}{\alpha_n(\si)}\frac{\alpha_2(\si)}{\alpha_n(\si)}\cdots\frac{\alpha_{n-1}(\si)}{\alpha_n(\si)} \times  n^{(n-1)/2}  \frac{\delta}{\alpha_1(\si)}\frac{\delta}{\alpha_2(\si)}\cdots\frac{\delta}{\alpha_{n-1}(\si)} \leq \kappa 2^n n^{(n-1)/2}  b^{-(n-1)}
\]
using \eqref{bt}. This contradicts the assumption that $\inf\limits_{\pi\in P_{n-1}}\mathcal{L}^{n-1}(\pi C) > 0$ since we can choose $\kappa$ arbitrarily small. \hfill $\square$

\section*{Acknowledgments}
SAB thanks the Carnegie Trust for financially supporting this work.  JMF was financially supported by a Leverhulme Trust Research Fellowship (RF-2016-500) and an EPSRC Standard Grant (EP/R015104/1).

\affiliationone{% in this example, two authors share an institution
   Stuart A. Burrell\\
   School of Mathematics and Statistics
   University of St Andrews, St Andrews, KY16 9SS, UK
   \email{sb235@st-andrews.ac.uk}}
   
\affiliationone{% in this example, two authors share an institution
   Jonathan M. Fraser\\
   School of Mathematics and Statistics
      University of St Andrews, St Andrews, KY16 9SS, UK
   \email{jmf32@st-andrews.ac.uk}}

\begin{thebibliography}{99}

% 1. Replace 9 by 99 if 10 or more references
% 2. Use "\and" between author names below
% 3. Format:
%    \bibitem{key}
%    {\bibname J. Smith, A. Kent \and D. I. Olive}, `Name of Paper', {\em
%    Journal Name } volume:start-end (year).

\bibitem{fraserbaker}
{\bibname S. Baker, J.M. Fraser \and \'A. M\'ath\'e}, \emph{Inhomogeneous self-similar sets with overlaps}, {Ergodic Th.  Dynam. Syst.}, 39(1):1-18 (2019).

\bibitem{bhr}
{\bibname B.B\'ar\'any, M. Hochman \and A. Rapaport}, {\emph{Hausdorff dimension of planar self affine sets and measures}}, {Invent. Math.}, 216:601-659 (2019). 

\bibitem{barn}
{\bibname M.F. Barnsley \and S. Demko}, \emph{Iterated function systems and the global construction of fractals}, {Proc. R. Soc. Lond. Ser. A}, 399:243-275 (1985).

\bibitem{burrell}
{\bibname S.A. Burrell}, \emph{On the dimension and measure of inhomogeneous attractors}, {Real Anal. Exchange}, 44(1):199-216 (2019).

\bibitem{1988kjf}
{\bibname K.J. Falconer}, \emph{The Hausdorff dimension of self-affine fractals}, {Math. Proc. Camb. Phil. Soc.}, 103:339-350 (1988).

\bibitem{Falconer}
{\bibname K.J. Falconer}, Fractal Geometry: Mathematical Foundations and Applications, {\em John Wiley \& Sons}, (1990).

\bibitem{hd2}
{\bibname K.J. Falconer}, \emph{The Hausdorff dimension of self-affine fractals II}, {Math. Proc. Camb. Phil. Soc.}, 111:169-179 (1992).

\bibitem{feng}
{\bibname D. Feng \and H. Hu}, \emph{Dimension theory of iterated function systems}, {Comm. Pure App. Math.}, 62(11):1435-1500 (2009)

\bibitem{fraserbox}
{\bibname J.M. Fraser}, \emph{On the packing dimension of box-like self-affine sets in the plane}, {Nonlinearity}, 15(1):77-97 (2012).

\bibitem{Fraser}
{\bibname J.M. Fraser}, \emph{Inhomogeneous self-similar sets and box dimensions}, {Studia Math.}, 213:133-156 (2012).

\bibitem{18}
{\bibname J.M. Fraser}, \emph{Inhomogeneous self-affine carpets}, {Indiana Univ. Math. J.}, 65:1547-1566 (2016).

\bibitem{fraserkemp}
{\bibname J.M. Fraser \and T. Kempton}, \emph{On the $L^q$ dimensions of measures on Heuter-Lalley type self-affine sets}, {Proc. Amer. Math. Soc.}, 146:161-173 (2018).

\bibitem{hata}
{\bibname M. Hata}, \emph{On some properties of set-dynamical systems}, {Proc. Japan Acad. Ser. A Math. Sci.}, 61:99-102 (1985).
 
\bibitem{27}
{\bibname I. Hueter \and S.P. Lalley}, \emph{Falconer's formula for the Hausdorff dimension of a self-affine set in $\R^2$},  {Ergodic Th. Dynam. Syst.}, 15(1):77-97 (1995).

\bibitem{jordanjurga}
{\bibname T. Jordan and N. Jurga}, \emph{Self affine sets with non-compactly supported random perturbations}, {Ann. Acad. Sci. Fen. Math.}, 39:771-785 (2014).

\bibitem{antti}
{\bibname A. K\"aenm\"aki \and J. Lehrb\"ack}, \emph{Measures with predetermined regularity and inhomogeneous self-similar sets}, {Ark. Mat.}, 116(4):929-956 (2018).

\bibitem{olssni}
{\bibname L. Olsen \and N. Snigireva}, \emph{$L^q$ spectra and R\'enyi dimensions of in-homogeneous self-similar measures}, {Nonlinearity}, 20:151-175 (2007).

\bibitem{sniphd}
{\bibname N. Snigireva}, Inhomogeneous self-similar sets and measures, {\em PhD Dissertation, University of St Andrews},  (2008).

\end{thebibliography}
\end{document}